\documentclass{amsart}
\usepackage{amssymb,amsmath}
\usepackage[mathscr]{euscript}

\newcounter{punct}

\def\punct{\refstepcounter{punct}{
\arabic{punct}.  }}

\newtheorem{theorem}{Theorem}
\newtheorem{proposition}[theorem]{Proposition}

\newtheorem{lemma}[theorem]{Lemma}

          \def\sm{\smallskip}

\begin{document}

\newcommand{\supp}{\mathop {\mathrm {supp}}\nolimits}
\newcommand{\rk}{\mathop {\mathrm {rk}}\nolimits}
\newcommand{\Aut}{\mathop {\mathrm {Aut}}\nolimits}
\newcommand{\Ob}{\mathop {\mathrm {Ob}}\nolimits}
\newcommand{\Out}{\mathop {\mathrm {Out}}\nolimits}
\renewcommand{\Re}{\mathop {\mathrm {Re}}\nolimits}
\newcommand{\Inn}{\mathop {\mathrm {Inn}}\nolimits}
\newcommand{\Char}{\mathop {\mathrm {Char}}\nolimits}
\newcommand{\Sp}{\mathop {\mathrm {Sp}}\nolimits}
\newcommand{\SOS}{\mathop {\mathrm {SO^*}}\nolimits}
\newcommand{\Ams}{\mathop {\mathrm {Ams}}\nolimits}
\newcommand{\Gms}{\mathop {\mathrm {Gms}}\nolimits}

\newcommand{\edge}{\mathop {\mathrm {edge}}\nolimits}

\def\0{\mathbf 0}

\def\ov{\overline}
\def\wh{\widehat}
\def\wt{\widetilde}
\def\pol{\twoheadrightarrow}

\renewcommand{\rk}{\mathop {\mathrm {rk}}\nolimits}
\renewcommand{\Aut}{\mathop {\mathrm {Aut}}\nolimits}
\renewcommand{\Re}{\mathop {\mathrm {Re}}\nolimits}
\renewcommand{\Im}{\mathop {\mathrm {Im}}\nolimits}
\newcommand{\sgn}{\mathop {\mathrm {sgn}}\nolimits}
\newcommand{\gr}{\mathop {\mathrm {gr}}\nolimits}

\def\bfa{\mathbf a}
\def\bfb{\mathbf b}
\def\bfc{\mathbf c}
\def\bfd{\mathbf d}
\def\bfe{\mathbf e}
\def\bff{\mathbf f}
\def\bfg{\mathbf g}
\def\bfh{\mathbf h}
\def\bfi{\mathbf i}
\def\bfj{\mathbf j}
\def\bfk{\mathbf k}
\def\bfl{\mathbf l}
\def\bfm{\mathbf m}
\def\bfn{\mathbf n}
\def\bfo{\mathbf o}
\def\bfp{\mathbf p}
\def\bfq{\mathbf q}
\def\bfr{\mathbf r}
\def\bfs{\mathbf s}
\def\bft{\mathbf t}
\def\bfu{\mathbf u}
\def\bfv{\mathbf v}
\def\bfw{\mathbf w}
\def\bfx{\mathbf x}
\def\bfy{\mathbf y}
\def\bfz{\mathbf z}

\def\bfA{\mathbf A}
\def\bfB{\mathbf B}
\def\bfC{\mathbf C}
\def\bfD{\mathbf D}
\def\bfE{\mathbf E}
\def\bfF{\mathbf F}
\def\bfG{\mathbf G}
\def\bfH{\mathbf H}
\def\bfI{\mathbf I}
\def\bfJ{\mathbf J}
\def\bfK{\mathbf K}
\def\bfL{\mathbf L}
\def\bfM{\mathbf M}
\def\bfN{\mathbf N}
\def\bfO{\mathbf O}
\def\bfP{\mathbf P}
\def\bfQ{\mathbf Q}
\def\bfR{\mathbf R}
\def\bfS{\mathbf S}
\def\bfT{\mathbf T}
\def\bfU{\mathbf U}
\def\bfV{\mathbf V}
\def\bfW{\mathbf W}
\def\bfX{\mathbf X}
\def\bfY{\mathbf Y}
\def\bfZ{\mathbf Z}

\def\frD{\mathfrak D}
\def\frQ{\mathfrak Q}
\def\frS{\mathfrak S}
\def\frT{\mathfrak T}
\def\frL{\mathfrak L}
\def\frG{\mathfrak G}
\def\frb{\mathfrak b}
\def\frg{\mathfrak g}
\def\frh{\mathfrak h}
\def\frf{\mathfrak f}
\def\frk{\mathfrak k}
\def\frl{\mathfrak l}
\def\frm{\mathfrak m}
\def\frn{\mathfrak n}
\def\fro{\mathfrak o}
\def\frp{\mathfrak p}
\def\frq{\mathfrak q}
\def\frr{\mathfrak r}
\def\frs{\mathfrak s}
\def\frt{\mathfrak t}
\def\fru{\mathfrak u}
\def\frv{\mathfrak v}
\def\frw{\mathfrak w}
\def\frx{\mathfrak x}
\def\fry{\mathfrak y}
\def\frz{\mathfrak z}

\def\bfw{\mathbf w}

\def\R {{\mathbb R }}
 \def\C {{\mathbb C }}
  \def\Z{{\mathbb Z}}
  \def\H{{\mathbb H}}
\def\K{{\mathbb K}}
\def\N{{\mathbb N}}
\def\Q{{\mathbb Q}}
\def\A{{\mathbb A}}

\def\T{\mathbb T}
\def\P{\mathbb P}
\def\SS{\mathbb S}

\def\G{\mathbb G}

\def\cD{\EuScript D}
\def\cL{\EuScript L}
\def\cK{\EuScript K}
\def\cM{\EuScript M}
\def\cN{\EuScript N}
\def\cP{\EuScript P}
\def\cT{\EuScript T}
\def\cQ{\EuScript Q}
\def\cR{\EuScript R}
\def\cW{\EuScript W}
\def\cY{\EuScript Y}
\def\cF{\EuScript F}
\def\cG{\EuScript G}
\def\cZ{\EuScript Z}
\def\cI{\EuScript I}
\def\cB{\EuScript B}
\def\cA{\EuScript A}
\def\cE{\EuScript E}
\def\cC{\EuScript C}
\def\cS{\EuScript S}

\def\bbA{\mathbb A}
\def\bbB{\mathbb B}
\def\bbD{\mathbb D}
\def\bbE{\mathbb E}
\def\bbF{\mathbb F}
\def\bbG{\mathbb G}
\def\bbI{\mathbb I}
\def\bbJ{\mathbb J}
\def\bbL{\mathbb L}
\def\bbM{\mathbb M}
\def\bbN{\mathbb N}
\def\bbO{\mathbb O}
\def\bbP{\mathbb P}
\def\bbQ{\mathbb Q}
\def\bbS{\mathbb S}
\def\bbT{\mathbb T}
\def\bbU{\mathbb U}
\def\bbV{\mathbb V}
\def\bbW{\mathbb W}
\def\bbX{\mathbb X}
\def\bbY{\mathbb Y}

\def\kappa{\varkappa}
\def\epsilon{\varepsilon}
\def\phi{\varphi}
\def\le{\leqslant}
\def\ge{\geqslant}

\def\B{\mathrm B}

\def\la{\langle}
\def\ra{\rangle}
\def\tri{\triangleright}

\def\lambdA{{\boldsymbol{\lambda}}}
\def\alphA{{\boldsymbol{\alpha}}}
\def\betA{{\boldsymbol{\beta}}}
\def\mU{{\boldsymbol{\mu}}}
\def\XI{{\boldsymbol{\Xi}}}

\def\const{\mathrm{const}}
\def\rem{\mathrm{rem}}
\def\even{\mathrm{even}}
\def\SO{\mathrm{SO}}
\def\SL{\mathrm{SL}}
\def\SU{\mathrm{SU}}
\def\PSL{\mathrm{PSL}}
\def\cont{\mathrm{cont}}
\def\O{\mathrm{O}}
\def\Th{\mathrm{Th}}

\def\U{\operatorname{U}}
\def\GL{\operatorname{GL}}
\def\Mat{\operatorname{Mat}}
\def\End{\operatorname{End}}
\def\Mor{\operatorname{Mor}}
\def\Aut{\operatorname{Aut}}
\def\inv{\operatorname{inv}}
\def\red{\operatorname{red}}
\def\Ind{\operatorname{Ind}}
\def\dom{\operatorname{dom}}
\def\im{\operatorname{im}}
\def\md{\operatorname{mod\,}}
\def\indef{\operatorname{indef}}
\def\Gr{\operatorname{Gr}}
\def\Pol{\operatorname{Pol}}
\def\Cut{\operatorname{Cut}}
\def\Add{\operatorname{Add}}
\def\ord{\operatorname{ord}}
\def\Replace{\operatorname{Replace}}
\def\Tr{\operatorname{Tr}}
\def\Homeo{\operatorname{Homeo}}
\def\Sep{\operatorname{Sep}}

\def\arr{\rightrightarrows}
\def\bs{\backslash}

\def\cH{\EuScript{H}}
\def\cO{\EuScript{O}}
\def\cQ{\EuScript{Q}}
\def\cL{\EuScript{L}}
\def\cX{\EuScript{X}}

\begin{center}
\large\bf
 
 On  algebras of double cosets of  symmetric groups
\\ 
with respect to  Young subgroups

 \medskip

\sc Yury A. Neretin%
	\footnote{The work is supported by the grant  FWF,  P31591.}
 
\end{center}

\bigskip

{\small 
We consider the subalgebra
 $\Delta$ in the group algebra of the symmetric group 
 $G=S_{n_1+\dots+n_\nu}$ consisting of all functions invariant
with respect to left and right shifts by elements of the Young subgroup 
 $H:=S_{n_1}\times \dots \times S_{n_\nu}$.
 We discuss structure constants of the algebra
 $\Delta$ and construct an algebra with continuous parameters $n_1$
 extrapolating algebras $\Delta$, it can be also can be rewritten as
 an asymptotic algebra as  $n_j\to\infty$ (for fixed $\nu$).
We show that there is a natural map from the Lie algebra of the group of pure braids
to 
 $\Delta$  
(and therefore this Lie algebra acts in spaces of multiplicities
of the quasiregular representation of the group  $G$ in functions
on  $G/H$).

}


\bigskip

{\bf\punct Introduction.}
Let 
 $G$ be a locally compact group, $K\subset G$ its compact subgroup.
 Consider the corresponding Hecke--Iwahori algebra,
i.e, the convolution algebra  
$\cM(G,K)$ of complex-valued compactly supported measures on $G$, which are
invariant with respect to left and right shifts by elements of   $K$ (see,
e.g., \cite{Bum}). We can regard such measures as measures on the space
of double cosets
$K\backslash G/K$.

It is easy to see that for any unitary representation of the group
  $G$ the algebra 
 $\cM(G,K)$ acts in the space of all  
$K$-invariant vectors. For this reason, Hecke-Iwahori algebras are a tool for 
investigation of
representations of groups. They are used in numerous contexts
by different reasons starting work by Gelfand  \cite{Gel},\
who introduced the notion of spherical subgroups and showed that the phenomenon
of sphericity is related to commutativity of an algebra
$\cM(G,K)$. Such algebras can be interesting objects by themselves and
can live by their own lives as the Hecke algebra, the affine Hecke algebra 
(both objects were introduced by Iwahori, see 
 \cite{Iwa}), the Yokonuma algebra  \cite{Yok}.

On the other hand, such algebras can be complicated objects
in cases that a priori seem simple. For instance, if
 $G=\SL(2,\R)$, $K=\SO(2)$ (or, more generally  $G$ is a real simple Lie group
 of rank 1 and $K$ is the maximal compact subgroup),
 then the product in $\cM(G,K)$ is determined
 by a hypergeometric kernel of type  $_2F_1$, see, e.g.,   \cite{Koo}, Sect.~7.
 In this case  
$\cM(G,K)$ is an interesting useful object, see Flensted-Jensen, Koornwinder
 \cite{FJ}, but a generalization of the mentioned hypergeometric formula to classical groups of rank $>1$ does not seem to be known
 (on some continuations of this topic, see 
 \cite{BG}, \cite{RV}). An example of a complicated object admitting nontrivial
 descriptions is the algebra of conjugacy classes of the symmetric group
 $S_n$, see Ivanov, Kerov
 \cite{IK}; on the algebra of conjugacy classes of  $\GL$ over a finite field,
 see M\'eliot  \cite{Mel}.

Olshanski \cite{Olsh} discovered that often a pass to an infinite limit
leads to  degenerations of algebras  $\cM(G_n,K_n)$ (as $\cM(\U(n),\O(n))$ or 
$\cM(S_{2n},S_n\times S_n)$);
on the set $K_\infty\backslash G_\infty/K_\infty$ itself  a
natural multiplication can arise. This, in its turn, creates a tool
for investigation of representations of infinite-dimensional groups.
The author (e.g., \cite{Ner}) observed that such multiplications 
are a mass phenomenon and admit explicit descriptions in a wide generality
(including descriptions of spaces $K_\infty\backslash G_\infty/K_\infty$).
The present note together with  \cite{Ner-tri} is a kind of an experiment --
an attempt of a descent from a relatively understandable infinite-dimensional
limit to prelimit objects. Clearly, our object is very complicated, we choose it
since in this case a limit picture seems the most transparent, see
 \cite{Nes}, \cite{Ner}, \S 8. 

\sm

{\bf\punct  Algebra of double cosets with respect to a Young subgroup. %
\label{ss:young}}
By $S_N$ we denote the group of all permutations of a set $\Omega$ with   $N$
elements. Let  
$$N=n_1+n_2+\dots + n_\nu.$$
Split $\Omega$ into a disjoint union  $\Omega=\coprod P_j$,
where a number  $\# P_j$ of elements of a set  $P_j$ is $n_j$.
Denote by
 $Y\{n_j\}$ the corresponding {\it Young subgroup}, i.e., the subgroup
in  $S_N$ consisting of permutations sending each set  $P_j$ to itself.
Obviously, $Y\{n_j\}\simeq \prod_j S_{n_j}$. To each double coset $ Y\{n_j\}\cdot g \cdot Y\{n_j\}$ we assign
an integer matrix 
$$a_{ij}=\# \bigl(g(P_i)\cap P_j\bigr). $$
It is easy to show (see, e.g., \cite{JK}, Sect. 1.3), that the matrix $\{a_{ij}\}$
satisfies conditions
\begin{equation}
\sum\nolimits_i a_{ij}=n_j, \quad \sum\nolimits_j a_{ij}=n_i, \quad a_{ij}\ge 0,
\label{eq:bist}
\end{equation}
and a double coset is uniquely determined by such matrix. 
We will use three notations for this double coset:

\sm

(1) $\xi\{a_{ij}\}$;

\sm

(2) $\xi\bigl(\sum a_{ij} E_{ij}\bigr)$ (where $E_{ij}$
is a matrix with unit at   $ij$-th place,
and all other matrix elements are 0);

\sm

(3) $\xi[A]$, where $A$ is an explicitly written matrix.

\sm

Obviously, the number of elements of a double coset
 $\xi\{a_{ij}\}$ is 
\begin{equation}
\mu\{a_{ij}\}=\frac{\prod (n_j!)^2}{\prod a_{ij}!}.
\label{eq:mu}
\end{equation}

{\sc Remark.}
The number of double cosets (i.e., the number of integer matrices
 $\{a_{ij}\}$ satisfying 
 \eqref{eq:bist}) is a complicated combinatorial function,
it was examined, in particular, in  
 \cite{Jon}, \cite{Dia}.
 \hfill $\boxtimes$
 
 \sm

Consider the group algebra
 $\C[S_N]$ and its element 
$$
\Pi:=\frac1{\prod n_j!}\sum_{h\in Y\{n_j\}} h.
$$ 
Clearly, $\Pi^2=\Pi$, and the subspace   
$$
\Delta\{n_j\}:=\Pi\, \C[S_N]\,\Pi\subset \Pi[S_N]
$$
coincides with the algebra of functions on $S_N$ invariant with respect to left and right shifts by elements of
 $Y\{n_j\}$.  A natural basis in this algebra consists of elements
$$
\Xi\{a_{ij}\}=\frac{1}{\mu\{a_{ij}\}} \sum_{g\in \xi\{a_{ij}\}} g.
$$
For basis elements
 $\Xi$ we use three types of notations
similarly to  notations (1)--(3)
for
 $\xi$.

\begin{proposition}
A product of basis elements of the algebra
 $\Delta\{n_j\}$
 is
\begin{multline}
\Xi\{a_{ij}\}\cdot \Xi\{b_{jk}\}=\\=
\sum_{\{c_{ik}\}}
\frac{\prod a_{ij}!\prod b_{jk}!}{\prod n_j!}
\Biggl(\sum_{\substack{\{t_{ijk}\}:\,t_{ijk}\ge 0,\\
\sum\limits_i t_{ijk}=b_{jk},\,\sum\limits_j t_{ijk}=c_{ik},\,
\sum\limits_k t_{ijk}=a_{ij}
}} \frac{1}{\prod_{ijk} t_{ijk}!}\Biggr)\, \Xi\{c_{ik}\}.
\label{eq:abc}
\end{multline}
\end{proposition}

{\sc Proof.} Let us evaluate the number of pairs
 $g\in \xi\{a_{ij}\}$, $h\in \xi\{b_{jk}\}$ such that
$hg\in \xi\{c_{ik}\}$. Set
$$
T_{ijk}:=g(P_i)\cap P_j\cap h^{-1} (P_k),\qquad t_{ijk}:=\#T_{ijk}.
$$
Fix numbers
 $t_{ijk}$. They satisfy  relations
\begin{align}
a_{ij}=\#\bigl(g(P_i)\cap P_j\bigr)=\#\bigl(\coprod_{k}T_{ijk}\bigr)
=\sum_{k} t_{ijk};
\label{eq:sum-a}
\\
b_{jk}=\#\bigl(h(P_j)\cap P_k\bigr)=\#\bigl(\coprod_{i}T_{ijk}\bigr)
=\sum_{i} t_{ijk};
\label{eq:sum-b}
\\
c_{ik}=\#\bigl(hg(P_i)\cap P_k\bigr)=\#\bigl(\coprod_{j}T_{ijk}\bigr)
=\sum_{j} t_{ijk}.
\label{eq:sum-c}
\end{align}

By the construction,
 $P_j=\coprod_{ik} T_{ijk}$, a collection of such partitions
 can be chosen by 
\begin{equation}
\prod_j \frac{n_j!}{\prod_{ik} t_{ijk}!}
\label{eq:prod1}
\end{equation}
ways. On the other hand,  
  $P_i=\coprod_{jk} g^{-1}T_{ijk} $, a collection of such partitions
 can be chosen by 
 \begin{equation*}
 \prod_i \frac{n_i!}{\prod_{jk} t_{ijk}!},
 \end{equation*}
 ways, this number coincides with
  \eqref{eq:prod1}. A collection of bijections 
 $g:g^{-1} T_{ijk}\to T_{ijk}$ can be chosen by $\prod_{ijk} t_{ijk}!$
 ways. Similarly, we choose sets
  $h T_{ijk}$ and bijections 
 $h:T_{ijk}\to T_{ijk}$. Thus, the desired number of  pairs
$g$, $h$ with fixed collection  $t_{ijk}$
  is
$$\Bigl( \frac{\prod_j n_j!}{\prod_{ijk} t_{ijk}!}\Bigr)^3 
\prod_{ijk} (t_{ijk}!)^2.
$$
To get a coefficient 
 \eqref{eq:abc}, we  must  sum these numbers over all collections
 $t_{ijk}$ satisfying
 \eqref{eq:sum-a}--\eqref{eq:sum-c} and divide by a normalization 
 factor  $\mu\{a_{ij}\}\, \mu\{b_{jk}\}$, see \eqref{eq:mu}.
\hfill $\square$.

\sm

{\sc Example: $\nu=2$.} Consider the case  $\nu=2$, $N=n_1+n_2$.
It is known that in this case the algebra
 $\Delta(n_1,n_2)$ is commutative and the pair  $S_N\supset S_{n_1}\times S_{n_2}$ is spherical, see, e.g.,
\cite{Cec}, \S 6.2,  spherical functions are expressed in terms
of Hahn hypergeometric orthogonal polynomials of type $_3F_2(1)$. 
The basis of  $\Delta(n_1,n_2)$ consists of functions 
$$
\phi_a=\Xi\begin{bmatrix}
n_1-a&a\\a&n_2-a
\end{bmatrix}, \qquad \text{where $a=0$, $1$, \dots, $\min(n_1,n_2)$.}
$$ 
Products have a form
 $\phi_a \phi_b=\sum s_{ab}^c \phi_c$, where
the structure constants are given by the formula
\begin{multline*}
s_{ab}^c=\frac{(a!\,b!)^2 (n_1-a)!\,(n_2-a)!\, (n_1-b)!\,(n_2-b)!}{n_1!\,n_2!}
\times\\\times
\!\!\!\!
\sum\limits_{\substack{\sigma,\tau:\\
0\le\sigma,\tau\le\min(a,b),\\
\sigma+\tau= a+b-c,
\\
n_1-a-b+\tau\ge 0,\\ n_2-a-b+\sigma\ge0
}}\!\!\!\!\!\!\!\!\!
\frac{1}{\sigma!\,\tau!\,(a-\sigma)! (b-\sigma)! \,(a-\tau)!\,(b-\tau)!\,
(n_1-a-b+\tau)!
(n_2-a-b+\sigma)!},
\end{multline*}
in \eqref{eq:abc}  we denote   $t_{121}=\sigma$, $t_{212}=\tau$.
Inequalities for 
$\sigma$, $\tau$ can be omitted if to assume that a factorial of a negative
integer is   $\infty$. Substituting  $\tau=a+b-c-\sigma$
we come to summation in  $\sigma$. The summands take the form 
$$
\frac{1}{\sigma!(a\!+\!b\!-\!c\!-\!\sigma)!(a\!-\!\sigma)! (b\!-\!\sigma)! 
(c\!-\!b\!+\!\sigma)!
(c\!-\!a\!+\!\sigma)!\,
(n_1-\!c\!-\!\sigma)!
(n_2-\!a\!-\!b\!+\!\sigma)!}.
$$
We transform factorials as
 $(p+\sigma)!=p!(p+1)_\sigma$,
$(q-\sigma)!=(-1)^{\sigma}q!/(-q)_{\sigma}$. If we perform this operation
for all factorials in the denominator    (this is possible if
 $a+b\le n_2$), then we come to 
\begin{multline*}
s_{a,b}^c=\frac{a!\, b!\,(n_1-a)!\,(n_2-a)!\, (n_1-b)!\,(n_2-b)!}
{n_1!\,n_2!\,(n_1-c)!\,(n_2-a-b)!\,(a+b-c)!\, (c-a)!\, (c-b)!}
\times\\\times
{}_4F_3\left[\begin{matrix} 
-a,-b,c-a-b,c-n_1\\
c-b+1,c-a+1,n_2-a-b+1
\end{matrix};\,1
 \right],
  \end{multline*}
where $_4F_3$ is a generalized hypergeometric function
 (in fact the sum is finite, notice that all parameters of the upper row
of $_4F_3$ are integers $\le 0$).
If $a+b> n_2$, then we pass to a summation in  $\sigma'=\sigma-a-b+n_2$
and again come to an expression of the form
 $_4F_3(1)$. 

In spite of simplicity of this case, the identity  
$\sum_\gamma s^\gamma_{ab} s^\delta_{\gamma c}=\sum_\gamma s^\delta_{a\gamma}s^\gamma_{bc} $ for structure constants providing the associativity
of the algebra $\Delta(n_1,n_2)$ is an identity of an unusual type for functions
 $_4F_3(1)$.
 \hfill $\boxtimes$

\sm

{\bf \punct The homomorphism of the Lie algebra of braid group.}
Recall that according Maltsev 
\cite{Mal1} any finitely generated nilpotent group without torsion
can be canonically embedded to a nilpotent Lie group as a cocompact lattice, 
and therefore any such discrete group
determines a Lie algebra. The construction admits an extension 
to residually nilpotent groups
 \cite{Mal2}, in particular it is valid for the group
 of pure braids. The corresponding Lie algebra  $\mathrm{Br}_\nu$
 was constructed by Kohno \cite{Koh} (see another way of construction 
 in \cite{Xi}), it has generators
  $L_{ij}$, where $i,j=1$, \dots, $\nu$ and $i\ne j$, such that $L_{ij}=L_{ji}$;
  relations ('infinitesimal braid relations')  are
\begin{align*}
[L_{ij},L_{jk}+L_{ki}]=0,\qquad \text{where $i$, $j$, $k$ are pairwise distinct;}
\\
[L_{ij},L_{kl}]=0,\qquad \text{where $i$, $j$, $k$, $l$ are pairwise distinct.}
\end{align*}
Note (see, e.g., \cite{Xi}) that elements  $L_{\nu 1}$, $L_{\nu 2}$, \dots, $L_{\nu (\nu-1)}$
generate a free Lie algebra  $\mathrm{Free}_{\nu-1}$,
it is an ideal in 
 $\mathrm{Br}_\nu$.
The algebra  $\mathrm{Br}_\nu$ is a semi-direct product
of the subalgebra 
$\mathrm{Br}_{\nu-1}$ generated by $L_{ij}$ with $i,j<\nu$, 
and the ideal  $\mathrm{Free}_{\nu-1}$. 

Kohno \cite{Koh1},\cite{Koh2} also showed that if the algebra
 $\mathrm{Br}_\nu$ acts in a finite dimensional linear space,
 then we have in this space a representation of the group of pure braids,
 this representation  
 arises from the Knizhnik--Zamolodchikov connection.
 The algebra  $\mathrm{Br}_\nu$ is graded by positive degrees 1, 2, 3, \dots 
 (degrees of generators are  1), this gives a possibility to define
 the corresponding infinite-dimensional Lie group as a certain
 'submanifold' in a completed universal enveloping algebra of $\mathrm{Br}_\nu$;
 there are different ways to define such completions.
 In particular, there is a completion, which acts in all finite-dimensional
 representations of $\mathrm{Br}_\nu$ and contains the group of pure braides, \cite{Ner-Mal}.

Many nontrivial actions of the Kohno algebra are known, see, e.g., 
\cite{Var}--\cite{TL}, \cite{Ner-Mal}, in particular it acts in various 
spaces of multiplicities of representations. By the Frobenius duality,
the next theorem shows that 
  $\mathrm{Br}_\nu$ and the corresponding 'Lie group' act in a natural way in spaces of multiplicities of the quasiregular representation of  $S_N$ in functions on $S_N/ Y\{n_i\}$.

Let $i\ne j$. Consider the transposition  $r_{ij}$,
exchanging an element of  $P_i$
and an element of  $P_j$. Set 
$$
\frr_{ij}:=\Pi\, r_{ij}\,\Pi=
\Xi\Bigl( \sum_{k:\,k\ne i,j} n_k E_{kk}+ (n_i-1)E_{ii}+(n_j-1)E_{jj}+E_{ij}+E_{ji}\Bigr).
$$

\begin{theorem}
Elements $\wt\frr_{ij}:=n_i n_j\, \frr_{ij}$
satisfy relations
\begin{align}
&[\wt \frr_{ij}, \wt \frr_{jk}+\wt \frr_{ik}]=0;
\label{eq:kon}
\\
&[\wt \frr_{ij}, \wt\frr_{kl}]=0, \quad \text{if $i$, $j$, $k$, $l$
are pairwise distinct.}
\label{eq:kon1}
\end{align}
\end{theorem}

{\sc Proof.} To verify \eqref{eq:kon} we without loss of generality
can assume 
$ijk=123$ and $\nu=3$,
\begin{multline*}
\frr_{23}\frr_{12}=\Xi
\begin{bmatrix}
n_1&0&0\\
0&n_2-1&1\\
0&1&n_3-1
\end{bmatrix}\,
\Xi\begin{bmatrix}
n_1-1&1&0\\
1&n_2-1&0\\
0&0&n_3-1
\end{bmatrix}=
\\
=\frac{n_2-1}{n_2}\cdot
\Xi\begin{bmatrix}
n_1-1&1&0\\
1&n_2-2&1\\
0&1&n_3-1
\end{bmatrix}
+
\frac{1}{n_2}\cdot
\Xi\begin{bmatrix}
n_1-1&0&1\\
1&n_2-1&0\\
0&1&n_3-1
\end{bmatrix}.
\end{multline*}
Therefore,
$$
[\frr_{23},\frr_{12}]=
\frac{1}{n_2}\cdot
\Xi\begin{bmatrix}
n_1-1&0&1\\
1&n_2-1&0\\
0&1&n_3-1
\end{bmatrix}-
\frac{1}{n_2}\cdot
\Xi\begin{bmatrix}
n_1-1&1&0\\
0&n_2-1&1\\
1&0&n_3-1
\end{bmatrix},
$$
and this implies
 \eqref{eq:kon}.

To establish \eqref{eq:kon1} it is sufficient to verify the identity for
 $ijkl=1234$ and $\nu=4$. We have
$$
\qquad\qquad
\frr_{12}\frr_{34}=\Xi\begin{bmatrix}
n_1-1&1&0&0\\
1&n_2-1&0&0\\
0&0&n_3-1&1\\
0&0&1&n_4-1
\end{bmatrix}=
\frr_{34}\frr_{12}.
\qquad\qquad\square
$$

\sm

{\bf\punct  Universalization of algebras $\Delta\{n_j\}$.}
Fix $\nu$. Denote by  $\cB$ the set of all integer
collections  
$\{a_{ij}\}$,   where $i$, $j\le \nu$ and $i\ne j$,
satisfying conditions  
\begin{equation*}
a_{ij}\ge 0, \qquad 
 \sum_{i:\,i\ne j} a_{ij}=\sum_{i:\,i\ne j} a_{j i}.
\end{equation*}
We define numbers
$$
a^*_{jj}:=\sum_{i:\,i\ne j} a_{j i}=\sum_{i:\,i\ne j} a_{ij}.
$$

We also introduce formal variables
 $\epsilon_1$, \dots, $\epsilon_\nu$.
Denote
$$
(\!(p,q;\epsilon_j)\!):=\prod_{m=p}^{q-1} (1-m\epsilon_j),
\qquad\text{in particular,}\quad (\!(p,p;\epsilon_j)\!):=1.
$$

For each collection
 $\{a_{ij}\}\in \cB$ we assign a formal basis element 
$
\XI\{a_{ij}\}
$.
Fix three such elements
 $\XI\{a_{ij}\}$, $\XI\{b_{jk}\}$, $\XI\{c_{ik}\}$. 
 Consider the set 
  $\cT=\cT(\{a_{ij}\},\{b_{jk}\},\{c_{ik}\})$
  of all collections  $\{t_{ijk}\}$ of non-negative
  integers, such that not all subscripts   $i$, $j$, $k$ coincide,
 satisfying the conditions 
\begin{align}
\sum\nolimits_{i} t_{ijk}&=b_{jk}\,\,\text{if}\,\, j\ne k;
\label{eq:bt}
\\
\sum\nolimits_{j} t_{ijk}&=c_{ik}\,\,\text{if}\,\, i\ne k;
\label{eq:ct}
\\
\sum\nolimits_{k} t_{ijk}&=a_{ij}\,\,\text{if}\,\, i\ne j,
\label{eq:at}
\end{align}
and also
\begin{equation}
a_{jj}^*+\sum_{k:k\ne j} t_{jjk}=b_{jj}^*+\sum_{i:i\ne j} t_{ijj}=
c_{jj}^*+\sum_{m:m\ne j} t_{jmj}=:t_{jjj}^*
\label{eq:t*}
\end{equation}
(we denote three equal values  \eqref{eq:t*} by $t_{jjj}^*$).
By these equalities,
\begin{align*}
\sum_{j,k:\,(j,k)\ne (i,i)} t_{ijk}=t^*_{iii};
\,\sum_{i,k:\,(i,k)\ne (j,)} t_{ijk}=t^*_{jjj};
\,\sum_{i,j:\,(i,j)\ne (k,k)} t_{ijk}=t^*_{kkk}.
\end{align*}

Note that sets
 $\cT(\{a_{ij}\},\{b_{jk}\},\{c_{ik}\})$
are bounded  in virtue of \eqref{eq:bt}--\eqref{eq:at} and $t_{ijk}\ge 0$,
therefore they are finite.
Moreover, for fixed
 $\{a_{ij}\}$, $\{b_{jk}\}$ the  set 
\begin{equation}
\bigcup\limits_{\{c_{ik}\}\in\cB} \cT(\{a_{ij}\},\{b_{jk}\},\{c_{ik}\}), 
\label{eq:fin}
\end{equation}
is finite.
Indeed, equalities
 \eqref{eq:bt}, \eqref{eq:at} are sufficient for boundedness.
Therefore for fixed   $\{a_{ij}\}$, $\{b_{jk}\}$
only finite number of sets
    $\cT(\{a_{ij}\},\{b_{jk}\},\{c_{ik}\})$ is non-empty.

Let us define the following values ({\it structure constants}):
\begin{multline*}
\frs_{\{a_{ij}\}, \{b_{jk}\}}^{ \{c_{ik}\}}=
\\=
\sum_\cT \frac{\prod\limits_{i,j:i\ne j} a_{ij}!\prod\limits_{j,k:j\ne k} b_{jk}!}
{\prod\limits_{\substack{i,j,k: \\\text{$i\ne j$ or $i\ne k$, or $j\ne k$}}}t_{ijk}!}
\prod_j 
 \epsilon_j^{a_{jj}^*+b_{jj}^*-t_{jjj}^*}
 \frac{(\!(a_{jj}^*,t_{jjj}^*;\epsilon_j)\!)
  (\!(b_{jj}^*,t_{jjj}^*;\epsilon_j )\!)}
{(\!(0,t_{jjj}^*;\epsilon_j )\!)}.
\end{multline*}

According the following lemma this expression is a well-defined
power series in
 $\epsilon_1$, \dots, $\epsilon_\nu$.

\begin{lemma}
\label{l:}
{\rm a)}
For any point $\in \cT(\{a_{ij}\},\{b_{jk}\},\{c_{ik}\})$ we have
  $a_{jj}^*+b_{jj}^*-t_{jjj}^*\ge 0$.
  
\sm  
  
 {\rm b)}  If $a_{jj}^*+b_{jj}^*-t_{jjj}^*= 0$
  for all $j$, then $a_{ij}+b_{ij}=c_{ij}$ for all pairs with  $i\ne j$.
\end{lemma}

{\sc Proof.}
$$
a^*_{jj}=\sum_{i:i\ne j} a_{ij} =\sum_{i:i\ne j}\sum_k t_{ijk}=t^*_{jjj}-
\sum_{k:k\ne j} t_{jjk}, \qquad b^*_{jj}=t^*_{jjj}-
\sum_{i:i\ne j} t_{ijj}.
$$
Therefore, 
$$
a_{jj}^*+b_{jj}^*-t_{jjj}^*=t^*_{jjj}-\sum_{k:k\ne j} t_{jjk}-
\sum_{i:i\ne j} t_{ijj}=\sum_{i,k:i\ne j, k\ne j} t_{ijk}.
$$
All summands in the right hand side are non-negative, therefore the sum is non-negative. This proves the statement a).

 If for all $j$ expressions
 $a_{jj}^*+b_{jj}^*-t_{jjj}^*$ are zero,  
 then only  nonzero $t_{\dots}$
can be $t_{ijj}$ or $t_{jjk}$. But in this case, $a_{ij}=\sum_k t_{ijk}=t_{ijj}$,
 $b_{ij}=\sum_{m}t_{mij}= t_{iij}$,
$c_{ij}=\sum_l t_{ilj}=t_{ijj}+t_{iij}$, and this proves
the statement  b).
\hfill $\square$

\sm

Consider two rings:

\sm
 
 --- First, denote by $\Q[[\{\epsilon_j\}]]$ the ring of
 formal series in
 $\epsilon_j$ with rational coefficients. 

\sm

--- Second,
let  $\Q\bigl[\{\epsilon_j\},\{(1-m\epsilon_j)^{-1}\}\bigr]$
be the ring of polynomial expressions in $\epsilon_j$, $\frac{1}{1-m\epsilon_j}$, where $j=1$, 
\dots, $\nu$ and $m=1$, $2$, \dots.

\sm

The structure constants
  $\frs^{\dots}_{\dots}$ can be regarded as elements of
  any of these rings.

\begin{theorem}
{\rm a)} The formula
$$
\XI\{a_{ij}\}\,\XI\{b_{jk}\}=\sum \frs_{\{a_{ij}\}, \{b_{jk}\}}^{ \{c_{ik}\}}
\XI\{c_{ik}\}
$$
determines a structure of an associative algebra over the ring
$\Q[[\{\epsilon_j\}]]$ and over the ring
$\Q\bigl[\{\epsilon_j\},\{(1-m\epsilon_j)^{-1}\}\bigr]$,
 denote this algebra by  $\Delta$.

\sm

{\rm b)} Fix $n_1$, \dots, $n_\nu$. 
Let 
\begin{equation}
a_{jj}^*, \,\,b_{jj}^*\le n_j.
\label{eq:le}
\end{equation}
Set
$$
\frs_{\{a_{ij}\}, \{b_{jk}\}}^{ \{c_{ik}\}}(\{n_j^{-1}\}):=
\frs_{\{a_{ij}\}, \{b_{jk}\}}^{ \{c_{ik}\}}
\Bigr|_{\epsilon_1=1/n_1,\, \dots,\,\epsilon_\nu= 1/n_\nu}.
$$
Then
$$
\frs_{\{a_{ij}\}, \{b_{jk}\}}^{ \{c_{ik}\}}(\{n_j^{-1}\})=0,\qquad \text{if
some
$c_{jj}^*>n_j$}.
$$
Moreover, the algebra
 $\Delta^\circ\{n_j\}$ with a basis $\XI^\circ \{a_{ij}\}$,
where $a_{jj}^*\le n_j$ and relations
$$
\XI^\circ\{a_{ij}\}\,\XI^\circ\{b_{jk}\}=
\sum_{\{c_{ik}\}} \frs_{\{a_{ij}\}, \{b_{jk}\}}^{ \{c_{ik}\}}(\{n_j^{-1}\})
\, \XI^\circ\{c_{ik}\}
$$
coincides with $\Delta\{n_j\}$.
\end{theorem}

{\sc Proof.} We write structure constants in formula
 \eqref{eq:abc} as
\begin{equation*} 
\prod_{i,j:i\ne j} a_{ij}! \,\prod_{j,k:j\ne k} b_{jk}!\, \sum 
\biggl(
\frac 1{\prod\limits_{\substack{i,j,k: \\\text{$i\ne j$ or $i\ne k$, or $j\ne k$}}}t_{ijk}!}
\cdot \prod_{j} \frac{a_{jj}!\,b_{jj}!}{n_j!\, t_{jjj}!} \biggr).
\end{equation*}

Next, we set $n_j-a_{jj}=:a_{jj}^*$, $1/n_j=\epsilon_j$
and transform the factor 
$\prod_j$ as
\begin{equation*}
\prod_j \frac{(n_j-a_{jj}^*)!\,(n_j-b_{jj}^*)!}{n_j! \,(n_j- t_{jjj}^*)!}.
\end{equation*}
Further,
\begin{multline*}
n_j!=n_j^{n_j}\cdot 1\cdot (1-\tfrac{1}{n_j})\cdot (1-\tfrac{2}{n_j})\dots
\cdot (1-\tfrac{n_j-1}{n_j})=\\=\epsilon_j^{-n_j} 1\cdot (1-\epsilon_j)(1-2\epsilon_j)
\dots (1-(n_j-1)\epsilon_j)=
\epsilon_j^{-n_j}(\!(0,t_{jjj}^*;\epsilon_j)\!) 
(\!(t_{jjj}^*, n_j;\epsilon_j)\!).
\end{multline*}
Similarly,
\begin{align*}
(n_j-t_{jjj}^*)!&=\epsilon_j^{-n_j+t_{jjj}^*}(\!(t_{jjj}^*, n_j;\epsilon_j)\!);
\\
(n_j-a_{jj}^*)!&=\epsilon_j^{-n_j+a_{jjj}^*}
(\!(a_{jj}^*, t_{jjj}^*;\epsilon_j)\!)
(\!(t_{jjj}^*, n_j;\epsilon_j)\!).
\end{align*}
Cancelling the fraction by $(\!(t_{jjj}^*, n_j;\epsilon_j)\!)^2$,
we come to the expression for  $\frs_{\dots}^{\dots}$.

\sm

Next, let all
 $a_{ll}^*$, $b_{ll}^*\le n_l$ and some  $t_{jjj}^*>n_j$.
 Then each of three factors of the expression
$$
 \frac{(\!(a_{jj}^*,t_{jjj}^*;\epsilon_j)\!)
  (\!(b_{jj}^*,t_{jjj}^*;\epsilon_j )\!)}
{(\!(0,t_{jjj}^*;\epsilon_j )\!)}
$$
contains $(1-n_j\epsilon_j)$, therefore  the substitution  $\epsilon_j=1/n_j$
gives 0. On the other hand, if some  $t_{kkk}^*\le n_k$,
then the factor  
$(1-n_k\epsilon_k)$ is absent in the numerator and in the denominator,
therefore we have no pole at
 $\epsilon_k=1/n_k$. Therefore, any summand, for which
 $t_{jjj}^*>n_j$ for some $j$, vanishes.

If for the same
    $a_{ll}^*$, $b_{ll}^*$ some  $c_{jj}^*>n_j$, then for all summands of
      $\frs_{\dots}^{\dots}$ we have $t_{jjj}^*\ge c_{jj}^* >n_j$, and therefore
      all summands are 0.

Thus we proved the statement 
 b). Next, consider the algebra over the ring 
$\Q\bigl[\{\epsilon_j\},\{(1-m\epsilon_j)^{-1}\}\bigr]$
defined in the statement 
a). Since the sets  \eqref{eq:fin} are finite,
both sides of the desired identity (associativity)
$$
(\XI\{a_{ij}\}\, \XI\{a_{jk}'\})\, \XI\{a_{kl}''\}  
= \XI\{a_{ij}\}\, (\XI\{a_{jk}'\}\, \XI\{a_{kl}''\}) 
$$
 expand in finite sums with rational coefficients
by the basis $\XI\{\cdot\}$. Substituting $\epsilon_j=1/n_j$ with
sufficiently large 
$n_j$ we get an identity in the algebra  $\Delta\{n_j\}$. Therefore, 
the desired identity for coefficients
$\in\Q\bigl[\{\epsilon_j\},\{(1-m\epsilon_j)^{-1}\}\bigr]$ also is fulfilled.
Therefore, expansions in series in  $\epsilon_j$ also coincide.
\hfill $\square$

\sm

{\sc Remark.} Substituting  $\epsilon_j\mapsto 1/n_j$
we get poles in some structure constants $\frs_{\dots}^{\dots}$ (as shows the example below).
Therefore we do not
get a well-defined {\it infinite-dimensional} algebra after this substitution
(at least for our normalization of basis elements).
But these poles do not appear under conditions \eqref{eq:le} of 
Theorem 4.b. Note also that after a substitution
 $\epsilon_j\mapsto \alpha_j$,  
where $\alpha_j\ne 1$, $\frac12$, $\frac13$, \dots for all 
$j$, we do not fall to poles of structure constants
and  get a well-defined infinite-dimensional
associative algebra.
\hfill $\boxtimes$

\sm

{\sc Example: $\nu=2$.} 
Let us return to the case $\nu=2$, discussed in Section 
\ref{ss:young}. The basis $\XI_a$ now is enumerated by
numbers  $a\ge 0$, the structure constants are defined by formulas
\begin{multline*}
\frs^c_{ab}=(a!\, b!)^2
\sum_{\substack{\sigma,\tau:\\
 0\le \sigma,\tau\le \min(a,b),\\
 \sigma+\tau=a+b-c}} 
 \frac{1}{\sigma!\,\tau!\,(a-\sigma)!\,(a-\tau)!\,(b-\sigma)!\,(b-\tau)!}
  \times\\ \times
 \epsilon_1^a \epsilon_2^a\cdot
 \frac{(\!(a,a+b-\tau; \epsilon_1)\!) \, (\!(b,a+b-\tau; \epsilon_1  )\!)}
 {(\!(0,a+b-\tau; \epsilon_1  )\!)}
\frac{(\!( a,a+b-\sigma; \epsilon_2 )\!) \, (\!(b,a+b-\sigma; \epsilon_2)\!)}
{(\!(0,a+b-\sigma; \epsilon_2   )\!)}.
\end{multline*}
Set $a=b$, $c=0$. Then the sum is reduced to a summand  with
 $\sigma=\tau=a$,
$$
\frs^0_{a,a}=\frac{(a!)^2 \cdot \epsilon_1^a \epsilon_2^a}{(\!(0,a;\epsilon_1)\!)\,(\!(0,a;\epsilon_2 )\!)}.
$$
If $a>n_1$ or $a>n_2$, then after the substitution  $\epsilon_1\to 1/n_1$,
$\epsilon_2\to 1/n_2$ to $\frs^0_{a,a}$ we get a pole.
\hfill $\boxtimes$ 

\sm

{\bf\punct The commutative limit and Poisson algebra.} Consider
the algebra  
$\Delta$ as an algebra over  $\Q[[\epsilon_1,\dots, \epsilon_\nu]]$.
Our algebra has a decreasing filtration in powers
of  
 $\epsilon$, \dots, $\epsilon_\nu$. By Lemma
 \ref{l:}, the corresponding graded algebra
   $\gr(\Delta)$ is commutative,
 the product is given by the formula
 $$
 \XI\{a_{ij}\} \,\XI\{b_{ij}\}=\XI\{a_{ij}+b_{ij}\}.
 $$
 I.e., we get a commutative algebra over
   $\Q$ with a basis $\XI\{a_{ij}\}$.
In fact this limit participates in considerations 
of  \cite{Nes}, \cite{Ner}, \S8.
 
 It is easy to write the next term of expansions  in powers
 of $\epsilon_j$, 
 \begin{multline*}
 \XI\{a_{ij}\} \,\XI\{b_{ij}\}=\XI\{a_{ij}+b_{ij}\}
 +\\
 +\sum_j \epsilon_j \sum_{\alpha,\gamma} a_{\alpha j} b_{j \beta}\, \XI\bigl(\sum_{i,k:i\ne k}
  (a_{ik}+b_{ik})E_{ik}+E_{\alpha \gamma}-E_{\alpha j}-E_{j\gamma}\bigr)
  +\sum_j O(\epsilon_j),
 \end{multline*}
 where $\sum_j O(\epsilon_j)$ denotes a formal series
 in monomials 
 $\prod_j \epsilon_j^{m_j}$ with $\sum m_j\ge 2$.
 
 Therefore,
 \begin{multline*}
 \XI\{a_{ij}\} \, \XI\{b_{ij}\} - \XI\{b_{ij}\} \, \XI\{a_{ij}\}   =
 \\=
 \sum_j \epsilon_j \sum_{\alpha,\gamma} 
 \bigl(a_{\alpha j} b_{j \beta}- a_{\beta j} b_{j \alpha} \bigr)\, \XI\bigl(\sum_{i,k:i\ne k}
  (a_{ik}+b_{ik})E_{ik}+E_{\alpha \gamma}-E_{\alpha j}-E_{j\gamma}\bigr)
  +\sum_j O(\epsilon_j).
 \end{multline*}
 In particular, assuming
  $\epsilon_1=\dots=\epsilon_\nu=\epsilon$, we get a Poisson bracket
 on the algebra $\gr (\Delta)$:
 $$
 \Bigl\{ \XI\{a_{ij}\} , \XI\{b_{ij}\}\Bigr\}:=
\biggl( \frac 1\epsilon\Bigl(\XI\{a_{ij}\} \, \XI\{b_{ij}\}-\XI\{b_{ij}\} \, \XI\{a_{ij}\}\Bigr)\biggr)\biggr|_{\epsilon=0},
 $$
 its explicit form is clear from the previous formula.

Math. Dept., University of Vienna;

Institute for Information Transmission Problems;

MechMath Dept., Moscow State University;

yurii.neretin(dog)univie.ac.at;

URL: http://mat.univie.ac.at/$\sim$neretin/

\end{document}